\documentclass[letterpaper, 10 pt, conference]{ieeeconf}
\IEEEoverridecommandlockouts
\usepackage{cite}
\usepackage{amsmath,amssymb,amsfonts}

\usepackage{mathtools}
\usepackage{graphicx}
\usepackage{mathbbol}
\usepackage{stfloats}
\usepackage{subfigure,color,balance,framed}
\usepackage{verbatim}
\usepackage{textcomp}
\usepackage{xcolor}
\usepackage{comment}
\usepackage{multirow}
\usepackage{makecell}
\usepackage{tabularx}

\usepackage{booktabs}
\usepackage{balance}
\usepackage{mathbbol}
\usepackage{algorithm}
\usepackage{algorithmicx}
\usepackage{algpseudocode}
\usepackage{mathrsfs}
\usepackage{threeparttable}

\newcommand{\method}[1]{\texttt{#1}}
\newcommand{\tr}{{{\mathsf T}}}

\usepackage[colorlinks=true,      
linkcolor=black,      
citecolor=black,      
filecolor=black,      
urlcolor=blue]{hyperref}
\newtheorem{definition}{Definition}
\newtheorem{theorem}{Theorem}
\newtheorem{proposition}{Proposition}
\newtheorem{remark}{Remark}

\def\BibTeX{{\rm B\kern-.05em{\sc i\kern-.025em b}\kern-.08em
    T\kern-.1667em\lower.7ex\hbox{E}\kern-.125emX}}

\title{\LARGE \bf Smoothing Mixed Traffic with Robust Data-driven Predictive Control for Connected and Autonomous Vehicles}

\author{Xu Shang$^{1}$, Jiawei Wang$^{2}$, and Yang Zheng$^{1}$
\thanks{The work of X. Shang and Y. Zheng is supported by NSF ECCS-2154650 and NSF CMMI-2320697.}
	\thanks{$^{1}$X. Shang and Y. Zheng are with the Department of Electrical and Computer Engineering, University of California San Diego, CA 92093, USA. (x3shang@ucsd.edu; zhengy@ucsd.edu), }%
	\thanks{$^{2}$J. Wang is with the Department of Civil and Environmental Engineering, University of Michigan, Ann Arbor, MI 48109, USA. (jiawe@umich.edu).}%
 }
\begin{document}
\maketitle

\begin{abstract}
The recently developed \method{DeeP-LCC} (Data-EnablEd Predictive Leading Cruise Control) method has shown promising performance for data-driven predictive control of Connected and Autonomous Vehicles (CAVs) in mixed traffic. However, its simplistic zero assumption of the future velocity errors for the head vehicle may pose safety concerns and limit its performance of smoothing traffic flow. In this paper, we propose a robust \method{DeeP-LCC} method to control CAVs in mixed traffic with enhanced safety performance. In particular, we first present a robust formulation that enforces a safety constraint for a range of potential velocity error trajectories, and then estimate all potential velocity errors based on the past data from the head vehicle. We also provide efficient computational approaches to solve the robust optimization for online predictive control. Nonlinear traffic simulations show~that our robust \method{DeeP-LCC} can provide better traffic efficiency and stronger safety performance while requiring less offline data.    
\end{abstract}


\section{Introduction}
In traffic flow, small perturbations of vehicle motion may propagate into large periodic speed fluctuations, leading to so-called stop-and-go traffic waves or phantom traffic jams~\cite{sugiyama2008traffic}. This phenomenon significantly lowers traffic efficiency and reduces driving safety. It has been widely demonstrated that connected and autonomous vehicles (CAVs) equipped with advanced control technologies, such as Cooperative Adaptive Cruise Control (CACC), have great potential to mitigate traffic jams \cite{milanes2013cooperative,li2017dynamical, zheng2015stability}. Yet, these technologies require a fully CAV environment,~and the near future will meet with a transition phase~of~mixed traffic where human-driven vehicles (HDVs) coexist with CAVs \cite{stern2018dissipation, zheng2020smoothing, orosz2016connected}. Thus, it is important to consider the behavior of HDVs when designing driving strategies for CAVs.

The control of CAVs in mixed traffic has indeed attracted increasing attention, and the existing methods are
generally categorized into model-based and model-free techniques. Model-based approaches typically use classical car-following models for HDVs, e.g., the Optimal Velocity Model (OVM)~\cite{bando1995dynamical}, to derive a parametric representation for mixed traffic. This parametric model is then utilized for CAV controller design, using methods such as optimal control~\cite{jin2016optimal, wang2021controllability}, $\mathcal{H}_{\infty}$ control \cite{ Mousavi2023MixedTraffic}, model predictive~control (MPC)~\cite{feng2021robust,zheng2016distributed}, and barrier methods \cite{zhao2023safety}. For these approaches, an accurate identification of the car-following models is non-trivial due to the complex and non-linear human driving behaviors. In contrast, model-free methods bypass system identification and directly design controllers for CAVs from data. For example, reinforcement learning \cite{wu2021flow} and adaptive dynamic programming \cite{huang2020learning} have been employed to learn wave-dampening CAV strategies. However, practical deployments of these methods are limited due to their computation burden and lack of interpretability and safety guarantees.  

Alternatively, data-driven predictive control methods that combine learning techniques with MPC have shown promising results for providing safe and optimal control of CAVs. In particular, the recent \method{DeeP-LCC} \cite{wang2023deep} exploits the Data-EnablEd Predictive Control (\method{DeePC}) \cite{coulson2019data,markovsky2021behavioral} technique for the Leading Cruise Control (LCC) \cite{wang2021leading} system in mixed traffic. This method directly utilizes the measured traffic data to design optimal control inputs for CAVs and explicitly incorporates input/output constraints in terms of limits on acceleration and car-following spacing. 
Large-scale numerical simulations~\cite{wang2023deep} and real-world experiments~\cite{wang2022implementation} have validated the capability of \method{DeeP-LCC} to smooth mixed traffic flow. However, the standard \method{DeeP-LCC} has an important zero velocity error assumption, i.e., the future velocity of the head vehicle remains the same as the equilibrium velocity of traffic flow. This assumption facilitates the online computation of \method{DeeP-LCC},   
but it will cause a mismatch between the real traffic behavior and its online prediction, which may compromise safety and 
control performance.

To address this issue, we develop a robust \method{DeeP-LCC} method to control CAVs in mixed traffic. Our key idea is to robustify \method{DeeP-LCC} by considering all potential velocity error trajectories and formulating a robust problem. 
We propose two methods for estimating velocity error trajectories and further present efficient computational approaches to solve the robust \method{DeeP-LCC} online via adapting standard robust optimization techniques \cite{bertsimas2011theory,lofberg2012automatic}. In particular, our main contributions include: 1) We propose a robust \method{DeeP-LCC} to handle the unknown velocity errors from the head vehicle. Our predictive controller will predict~a~series of future outputs based on the disturbance set and requires all of them to satisfy the safety constraint, thus providing enhanced safety performance. 2) We introduce two disturbance estimation methods, the constant velocity model and the constant acceleration model, based on the past disturbance data of the head vehicle. Our methods are able to provide good estimations of the future velocity errors and improve the control performance of robust \method{DeeP-LCC}. 
3) We further provide efficient computational approaches for solving the robust optimization problem. We analyze and compare the complexity of two different solving methods from the robust optimization literature \cite{bertsimas2011theory, lofberg2012automatic} and further provide a down-sampling method, adapted from \cite{huang2021decentralized}, to further decrease the computational complexity. Numerical experiments validate the enhanced performance of the robust \method{DeeP-LCC} in reducing fuel consumption and improving driving safety while requiring less pre-collected data. 
For~example, our robust \method{DeeP-LCC} only results in 4 and 0 emergencies out of 100 safety tests using small and large offline datasets,~respectively;  however, these numbers for \method{DeeP-LCC} \cite{wang2023deep} are 66 and 51 (which are unacceptably large).

The rest of the paper is organized as follows. Section \ref{sec:CF-LCC_DeeP-LCC} reviews the background on mixed traffic and \method{DeeP-LCC} for CF-LCC. Section \ref{sec:robust_DeeP-LCC} presents our robust \method{DeeP-LCC}. The disturbance set estimation methods and efficient computations are discussed in Section \ref{sec:Dis_Compute}. Section \ref{sec:Sim_Result} demonstrates our numerical results. We conclude the paper in Section \ref{sec:conclusion}. 

\section{Data-driven Predictive Control in CF-LCC}
\label{sec:CF-LCC_DeeP-LCC}
\begin{figure}[t]
    \centering
    \includegraphics[width=0.48\textwidth]{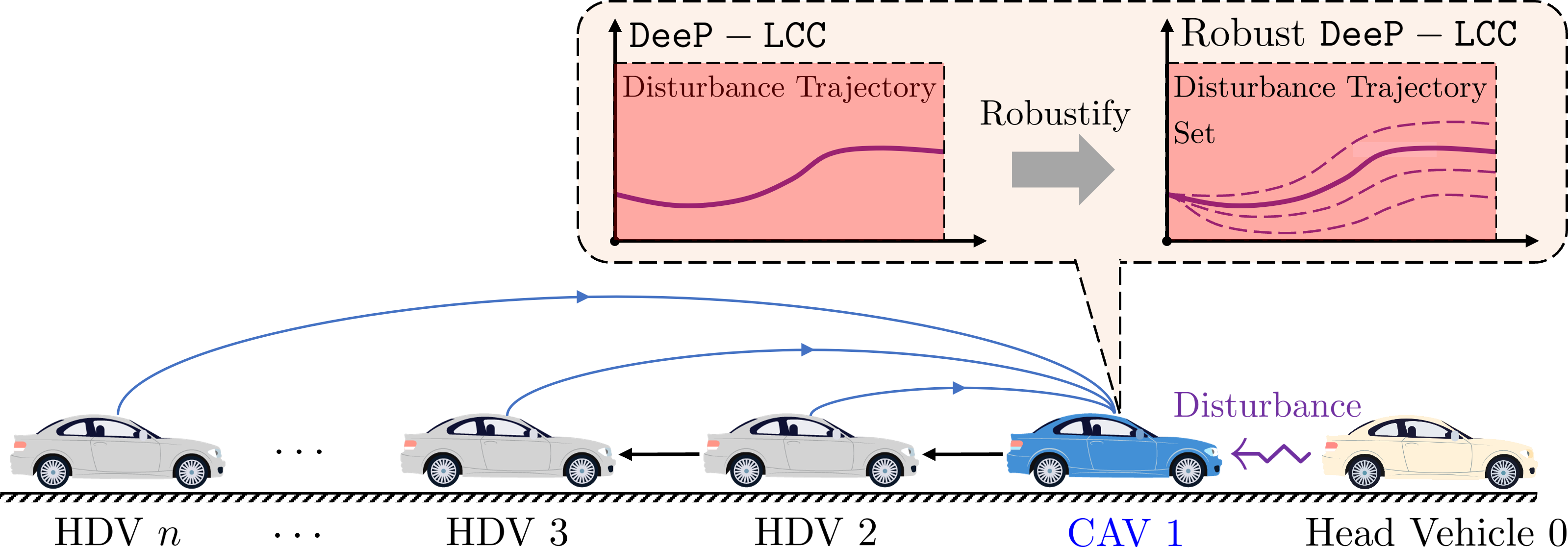}
    \caption{Schematic of CF-LCC system. 
    Original \method{DeeP-LCC} assumes one single trajectory for future disturbance, while the proposed robust \method{DeeP-LCC} explicitly addresses an estimated set of future disturbances.} 
    \label{fig:sys-CF-LCC}
    \vspace{-2mm}
\end{figure}

In this section, we briefly review the \method{DeeP-LCC}~\cite{wang2023deep} for a Car-Following LCC (CF-LCC) system~\cite{wang2021leading}. 
As shown in Fig. \ref{fig:sys-CF-LCC}, the CF-LCC consists of one CAV, indexed as $1$, and $n-1$ HDVs, indexed as $2,\ldots, n$ from front to end. All these vehicles follow a head vehicle, indexed as $0$, which is immediately ahead of the CAV. Such a CF-LCC system can be considered the smallest unit for general cascading mixed traffic~systems~\cite{wang2021leading}. Our robust \method{DeeP-LCC} can be extended to general mixed traffic systems; the details will be discussed in an extended report.

\subsection{Input/Output of CF-LCC system}

For the $i$-th vehicle at time $t$, we denote its position, velocity and acceleration as $p_i(t)$, $v_i(t)$ and $a_i(t)$, $i=1, \ldots, n$, respectively. We define the spacing between vehicle $i$ and its preceding vehicle as $s_i(t) = p_{i-1}(t)-p_i(t)$ and their relative velocity as $\dot{s}_i(t) = v_{i-1}(t)-v_i(t)$. In an equilibrium state, each vehicle moves at the same velocity $v^*$  with an equilibrium spacing $s_i^*$ that may vary from different vehicles. 

In \method{DeeP-LCC}, we consider the error state of the traffic system. In particular, the velocity error and spacing error for each vehicle are defined as
$
\tilde{v}_i(t) = v_i(t) - v^*,  \tilde{s}_i(t) = s_i(t) - s_i^*
$. 
Then, we form the state $x \in \mathbb{R}^{2n}$ of the CF-LCC system by lumping the error states of all the vehicles
\begin{equation*} \label{eq:sys-state}
x(t) = [\tilde{s}_1(t), \tilde{v}_1(t),\tilde{s}_2(t), \tilde{v}_2(t),\ldots,\tilde{s}_n(t), \tilde{v}_n(t)]^\tr.
\end{equation*}
The spacing errors of HDVs are not directly measurable, since it is non-trivial to get the equilibrium spacing $s_i^*$ for HDVs due to the unknown car-following behaviors. By contrast, the equilibrium velocity $v^*$ can be estimated from the past velocity trajectory of the leading vehicle. Accordingly, the system output is formed by the velocity errors of all vehicles and the spacing error of the CAV only, defined as
\[
    y(t) = [\tilde{v}_1(t), \tilde{v}_2(t),\ldots, \tilde{v}_n(t), \tilde{s}_1(t)]^\tr\in \mathbb{R}^{n+1}.
\]

The input $u(t) \in \mathbb{R}$ of the system is defined as the acceleration of the CAV, as widely used in \cite{orosz2016connected,zheng2020smoothing}. Finally, the velocity error of the head vehicle $0$ is regarded as an external disturbance signal $\epsilon = \tilde{v}_0(t) = v_0(t)-v^* \in \mathbb{R}$, and its past trajectory can be recorded, but its future trajectory is in general unknown.
Based on the definitions of the system state, input, and output, after linearization and discretization,  a state-space model of the CF-LCC system is in the form of
\begin{equation} \label{eq:Model-CF-LCC}
\left\{
\begin{aligned}
x(k+1) & = A x(k) + B u(k) + H \epsilon(k), \\
y(k) &= C x(k),
\end{aligned}
\right.
\end{equation}
where $k$ denotes the discrete time step. The details of the matrices $A, B, C, H$ can be found in~\cite[Section II-C]{wang2023deep}. 

Note that the parametric model \eqref{eq:Model-CF-LCC} is non-trivial to accurately obtain due to the unknown HDVs' behavior~(all different models, such as OVM, will lead to a system in the same form \eqref{eq:Model-CF-LCC}; see \cite{wang2023deep,zheng2020smoothing,wang2021leading,orosz2016connected} for details).~To address this issue, the recently proposed \method{DeeP-LCC} method directly uses the input/output trajectories for behavior prediction and controller design, thus bypassing the system identification process that is common in model-based methods.  

\subsection{Data-Driven Representation of System Behavior}

\method{DeeP-LCC} is an adaption of the standard \method{DeePC}~\cite{coulson2019data} for mixed traffic control. It starts by forming a data-driven representation of the system with rich enough pre-collected offline data and employs it as a predictor to predict the dynamical behavior of CF-LCC \eqref{eq:Model-CF-LCC}. We recall a persistent excitation \cite{willems2005note} for offline data collection.
\begin{definition}[Persistently Exciting]
    The sequence of signal $\omega = \textrm{col}(\omega(1),\omega(2), \ldots, \omega(T))$ with length $T$ ($T \in \mathbb{N}$) is persistently exciting of order $L$ ($L < T$) if its associated Hankel matrix with depth $L$ has full row rank:
    \[\mathcal{H}_L(\omega) = \begin{bmatrix}
        \omega(1) & \omega(2) & \cdots & \omega(T-L+1) \\
        \omega(2) & \omega(3) & \cdots & \omega(T-L+2) \\
        \vdots    & \vdots    & \ddots & \vdots \\
        \omega(L) & \omega(L+1) & \cdots & \omega(T)
    \end{bmatrix}.\]
\end{definition}

We begin with collecting an input/output trajectory of length $T$ for the CF-LCC system offline:
\[
\begin{aligned}
u^\textnormal{d} &= \textrm{col}(u^\textnormal{d}(1), u^\textnormal{d}(2),\ldots,u^\textnormal{d}(T))\in \mathbb{R}^{T}, \\
\epsilon^\textnormal{d} &= \textrm{col}(\epsilon^\textnormal{d}(1), \epsilon^\textnormal{d}(2),\ldots,\epsilon^\textnormal{d}(T))\in \mathbb{R}^{T},\\
y^\textnormal{d} &= \textrm{col}(y^\textnormal{d}(1), y^\textnormal{d}(2),\ldots,y^\textnormal{d}(T))\in \mathbb{R}^{(n+1)T}.
\end{aligned}
\] 
We then use the offline collected data to form a Hankel matrix of order $L$, which is partitioned as follows 
\begin{equation} \label{eq:Hankel-CF-LCC}
\begin{bmatrix}
    U_{\textnormal{P}} \\
    U_{\textnormal{F}} 
\end{bmatrix} \!:= \!\mathcal{H}_L(u^\textnormal{d}),\; 
\begin{bmatrix}
    E_{\textnormal{P}} \\
    E_{\textnormal{F}} 
\end{bmatrix}\! := \!\mathcal{H}_L(\epsilon^\textnormal{d}),\; 
\begin{bmatrix}
    Y_{\textnormal{P}} \\
    Y_{\textnormal{F}} 
\end{bmatrix} \!:= \!\mathcal{H}_L(y^\textnormal{d}),
\end{equation}
where $U_{\textrm{P}}$ and $U_{\textrm{F}}$ containts the first $T_{\textrm{ini}}$ rows and the last $N$ rows of $\mathcal{H}_L(u^{\textrm{d}})$, respectively (similarly for $E_\textrm{P}$ and $E_\textrm{F}$, $Y_\textrm{P}$ and $Y_\textrm{F}$).
The Hankel matrices \eqref{eq:Hankel-CF-LCC} can be used to construct the online behavior predictor for predictive control. Note that the CF-LCC system in \eqref{eq:Model-CF-LCC} is controllable; see a detailed proof in \cite{wang2021leading}.  
Then, we have the following result.
\begin{proposition}[\!\!{\cite[Proposition 2]{wang2023deep}}] \label{proposition:data representation}
At time step $k$, we collect the most recent past input sequence $u_{\textnormal{ini}}$ with length $T_{\textnormal{ini}}$, and let the future input sequence $u$ with length $N$ as 
\begin{align*}
 u_{\textnormal{ini}} &= \textrm{col}(u(k-T_\textnormal{ini}), u(k-T_\textnormal{ini}+1),\ldots,u(k-1)), \\
 u &= \textrm{col}(u(k), u(k+1), \ldots , u(k+N-1)).
\end{align*}
The notations $\epsilon_{\textnormal{ini}}$, $\epsilon$, $y_\textnormal{ini}$ and $y$ are denoted similarly.  
%
If the input trajectory $u^\textnormal{d}$ is persistently exciting of order $L+2n$ (where $L = T_\textnormal{ini} + N$), then 
the sequence $\textrm{col}(u_{\textnormal{ini}}, \epsilon_{\textnormal{ini}},$ $ y_{\textnormal{ini}},u,\epsilon,y)$ is a valid trajectory with length $L$ of \eqref{eq:Model-CF-LCC}  if and only if there exists a vector $g \in \mathbb{R}^{T-L+1}$ such that 
\begin{equation} \label{eq:CF-LCC-local-representation}
\begin{bmatrix}
U_{\textnormal{P}}\\
E_{\textnormal{P}}\\
Y_{\textnormal{P}}\\
U_{\textnormal{F}}\\
E_{\textnormal{F}}\\
Y_{\textnormal{F}}\\
\end{bmatrix} g
= 
\begin{bmatrix}
u_{\textnormal{ini}}\\
\epsilon_{\textnormal{ini}}\\
y_{\textnormal{ini}}\\
u\\
\epsilon\\
y\\
\end{bmatrix}.
\end{equation}
If $T_\textnormal{ini} \ge 2n$, then $y$ is unique for any $(u_{\textnormal{ini}}, y_{\textnormal{ini}}, u, \epsilon)$.
\end{proposition}
\vspace{1mm}

This proposition establishes a data-driven representation \eqref{eq:CF-LCC-local-representation} for the CF-LCC system: all valid trajectories can be~constructed by a linear combination of rich enough pre-collected trajectories. Thus, we can predict the future output $y$ using trajectories $(u^\textrm{d}, \epsilon^\textrm{d}, y^\textrm{d})$, given the future input $u$, disturbance $\epsilon$ and initial condition $(u_\textnormal{ini}, \epsilon_\textnormal{ini},y_\textnormal{ini})$.

\subsection{\method{DeeP-LCC} Formulation}
Using the data-driven representation \eqref{eq:CF-LCC-local-representation}, the \method{DeeP-LCC} in \cite{wang2023deep} solves an optimization problem at each time step:
\begin{subequations}\label{eqn:DeeP-LCC}
\begin{align}
\min_{g, \sigma_{y}, u, \epsilon, y}  \quad & V(u, y) + \lambda_{g} ||g||_2^2 + \lambda_{y} ||\sigma_{y}||_2^2 \label{eqn:objFun} \\
\textrm{subject~to} \quad & 
\begin{bmatrix}
U_{\textnormal{P}}\\
E_{\textnormal{P}}\\
Y_{\textnormal{P}}\\
U_{\textnormal{F}}\\
E_{\textnormal{F}}\\
Y_{\textnormal{F}}\\
\end{bmatrix} g
= 
\begin{bmatrix}
u_{\textnormal{ini}}\\
\epsilon_{\textnormal{ini}}\\
y_{\textnormal{ini}}\\
u\\
\epsilon\\
y\\
\end{bmatrix} + 
\begin{bmatrix}
0\\
0\\
\sigma_{y}\\
0\\
0\\
0\\
\end{bmatrix}\label{eqn:equality},\\
& \tilde{s}_{\min} \le G_1 y\le \Tilde{s}_{\max}, \label{eqn:safety}\\
& u_{\min} \le u \le u_{\max}, \label{eqn:inputlimit}\\
& \epsilon = \epsilon_{\textrm{est}} \label{eqn:estimator},
\end{align}
\end{subequations}
where $G_1 = I_N \otimes \begin{bmatrix}\mathbb{0}_{1\times n},\; 1\end{bmatrix}$ selects the spacing error of the CAV from the output, $[\tilde{s}_{\min}, \tilde{s}_{\max}]$ is the safe spacing error range of CAV, $[u_{\min}, u_{\max}]$ is the physical limitation of the acceleration and $\epsilon_{\textrm{est}}$ is the estimation of the future velocity errors of the head vehicle $0$. 

For the cost function \eqref{eqn:objFun}, $V(u, y)$ penalizes the output deviation from equilibrium states and the energy of the input:
\[
V(u, y) = ||u||_{R}^2 + ||y||_{Q}^2,
\] with $R \in \mathbb{S}_{+}^{N\times N}$ and $Q \in \mathbb{S}_+^{N(n+1)\times N(n+1)}$. 
There are two regularization terms $\|g\|_2^2$ and $\|\sigma_y \|_2^2$ in the cost function with weight coefficients $\lambda_g,\lambda_y$. Also, a slacking variable $\sigma_y$ is added to the data-driven representation~\eqref{eqn:equality}. Note that the original data-driven behavior representation~\eqref{eq:CF-LCC-local-representation} is only applicable to linear systems with noise-free data. The regularization herein is commonly used for nonlinear systems with stochastic noises, and we refer interested readers to~\cite{wang2023deep,coulson2019data} for detailed discussions. 

\begin{remark}[Robustification]
The \method{DeeP-LCC} \eqref{eqn:DeeP-LCC} requires an estimated sequence $\epsilon_{\textrm{est}}$ for the future disturbance (i.e., velocity errors) of the head vehicle.  In the standard~\method{DeeP-LCC}~\cite{wang2023deep}, it is assumed that \textit{the estimated future velocity error is zero}, which was justified by the assumption that vehicle $0$ always tries to maintain its equilibrium state. However, this assumption hardly stands since in real-world traffic, strong oscillations may happen, particularly during the occurrence of traffic waves. An inaccurate estimation of future velocity errors could cause a mismatch between the prediction and the real traffic behavior, which may not only degrade the control performance but also pose safety concerns (e.g., collision). In this paper, we will incorporate a valid set for disturbance estimation (see Fig. \ref{fig:sys-CF-LCC} for illustration) and establish a robust \method{DeeP-LCC}, as well as its tractable computations. 
    \hfill $\square$
\end{remark}

\section{Tractable Robust \method{DeeP-LCC} Formulation}
\label{sec:robust_DeeP-LCC}
In this section, we present a new framework of robust \method{DeeP-LCC} to control the CAV in the CF-LCC system which can properly address unknown future velocity errors, leading to enhanced performance and safety. 

\subsection{Robust \method{DeeP-LCC} Formulation}
As shown in Fig. \ref{fig:sys-CF-LCC}, instead of estimating one single disturbance trajectory $\epsilon = \epsilon_{\textrm{est}}$ in \method{DeeP-LCC}, we introduce a disturbance set $\mathcal{W}$ as the estimation, i.e., $\epsilon \in \mathcal{W}$, which, by valid design (see details in Section \ref{sec:Dis_Compute}), will contain the real trajectory with a much higher possibility. 

Our key idea is to plan over the worst trajectory~in~$\mathcal{W}$ for predictive control, 
leading to a robust optimization problem 
\begin{equation}\label{eqn:rDeeP}
\begin{aligned}
\min_{g, \sigma_{y}, u, y} \ \max_{\epsilon \in \mathcal{W}} \quad & V(u, y) + \lambda_{g} ||g||_2^2 + \lambda_{y} ||\sigma_{y}||_2^2  \\
\textrm{subject~to} \quad & \eqref{eqn:equality}, \ \eqref{eqn:safety}, \ \eqref{eqn:inputlimit}.
\end{aligned}
\end{equation}
Compared with the original formulation in \eqref{eqn:DeeP-LCC}, the robust formulation \eqref{eqn:rDeeP} promises to provide better control performance and higher safety guarantees since the gap between online prediction and real implementation has been reduced. 
As a trade-off, note that the complexity of the optimization problem is increased, and we also need to estimate $\mathcal{W}$ properly. Both issues will be discussed in the sections below. 

For implementation, the optimization problem \eqref{eqn:rDeeP} is solved in a receding horizon manner at each time step $k$ and $\mathcal{W}$ is re-estimated iteratively based on the updated velocity errors of the head vehicle (see Section \ref{sec:Dis_Compute}). Algorithm \ref{Alg:rDeeP-LCC} lists the overall procedure of robust \method{DeeP-LCC}.

\begin{algorithm}[t] 
	\caption{Robust \method{DeeP-LCC}}
	\label{Alg:rDeeP-LCC}
	\begin{algorithmic}[1]
		\Require
		Pre-collected offline data $(u^{\textrm{d}},\epsilon^{\textrm{d}},y^{\textrm{d}})$, initial time step $k_0$, terminal time step $k_f$;
		\State Construct data Hankel matrices \eqref{eq:Hankel-CF-LCC} for input, disturbance, and output as $U_{\textrm{P}}, U_{\textrm{F}}, E_{\textrm{P}}, E_{\textrm{F}}, Y_{\textrm{P}}, Y_{\textrm{F}}$;
		\State Initialize the most recent past traffic data $(u_{\textrm{ini}},\epsilon_{\textrm{ini}},y_{\textrm{ini}})$ before the initial time $k_0$;
		\While{$k_0 \leq k \leq k_f$}
            \State Estimate $\mathcal{W}$ from $\epsilon_{\textrm{ini}}$;
		\State Solve~\eqref{eqn:rDeeP} for optimal future control input sequence $u^*=\textrm{col}(u^*(k),u^*(k+1),\ldots,u^*(k+N-1))$;
		\State Apply the input $u(k) \leftarrow u^*(k)$ to the CAV;
		\State $k \leftarrow k+1$;
            \State Update past traffic data $(u_{ \textrm{ini}},\epsilon_{\textrm{ini}},y_{\textrm{ini}})$;
		\EndWhile
	\end{algorithmic}
\end{algorithm}

\subsection{Reformulations of the min-max optimization}

The min-max optimization problem \eqref{eqn:rDeeP} is solved at each iteration of Algorithm \ref{Alg:rDeeP-LCC}, but standard solvers are not applicable with respect to its current form. We proceed to present a sequence of reformulation (and relaxations) for \eqref{eqn:rDeeP}, which further allows for efficient computations in Section \ref{sec:Dis_Compute}. 

We first eliminate the equality constraint by expressing $g$ and $y$ in terms of $u$, $\sigma_y$ and $\epsilon$ to :
\begin{subequations}\label{eqn:g_rep}
\begin{align}
g &=
H_\textrm{p}^{\dag}b \!+\!H_\textrm{p}^{\bot}z \label{eqn:g_rep-a},\\
y &= Y_{\textnormal{F}} g = Y_{\textnormal{F}} H_\textrm{p}^\dag b + Y_{\textnormal{F}}H_\textrm{p}^{\bot}z, \label{eqn:g_rep-b}
\end{align}
\end{subequations}
where $H_\textrm{p} = \textrm{col}(U_{\textnormal{P}}, E_{\textnormal{P}}, Y_{\textnormal{P}}, U_{\textnormal{F}}, E_{\textnormal{F}})$ with $H_\textrm{p}^\dag$ denoting its pseudo-inverse, $H_\textrm{p}^\bot = I-H_\textrm{p}^\dag H_\textrm{p}$,  $b = \textrm{col}(u_{\textnormal{ini}}, \epsilon_{\textnormal{ini}}, y_{\textnormal{ini}}+\sigma_{y}, u, \epsilon)$, and $z \in \mathbb{R}^{T-L+1}$. For simplicity, we set $z = 0$ in the following derivation, which decreases the complexity of the optimization problem but also reduces the feasible set. From the simulations in Section \ref{sec:Sim_Result}, we note that this simplification already provides satisfactory control performance. Then, the min-max robust problem \eqref{eqn:rDeeP} 
becomes: 
\begin{subequations}
\label{eqn:robustCount}
\begin{align}
\min_{u, \sigma_{y}} \ \max_{\epsilon \in \mathcal{W}} \quad & x^\tr M x + d^\tr x + c_0 \label{eqn:Obj-robust}\\
\textrm{subject~to} \quad  
& \tilde{s}_{\min} \leq P_1 x + c_1 \leq \tilde{s}_{\max} \label{eqn:estimator-robust}, \\
& u_{\min} \le P_2 x \le u_{\max}, \label{eqn:inputlimit-robust}
\end{align}
\end{subequations}
where $x=\textrm{col}(u, \sigma_y, \epsilon) $ denotes the decision variable\footnote{With slight abuse of notations, we use $x$ to denote the decision variable in robust optimization.}, and {$M,d,c_0, P_1, P_2,c_1$ only depend on problem data (their explicit forms are provided in our numerical implementation).} 

Without loss of generality, we eliminate the constant~$c_0$. We finally consider $\epsilon$ as an uncertainty parameter,~and~transform problem \eqref{eqn:robustCount} into its epi-graph form 
\begin{subequations}
\label{eqn:epi_robustCount}
\begin{align}
\min_{x, t} & \quad t \nonumber \\
\textrm{subject~to}& \quad  x^\tr M x + d^\tr x \leq t, \quad \forall \epsilon \in \mathcal{W}, \label{eqn:obj_rb}\\ 
& \quad \tilde{s}_{\min} \le P_1 x + c_1\le \tilde{s}_{\max}, \quad \forall \epsilon \in \mathcal{W}, \label{eqn:safe_rb}\\
& \quad u_{\min} \le P_2 x \le u_{\max}.  
\end{align}
\end{subequations}
    Compared with \eqref{eqn:robustCount}, the formulation \eqref{eqn:epi_robustCount} requires its feasible solutions to satisfy the safety constraint for any $\epsilon$. This design indicates that the predictive controller needs to ensure safe constraints for all disturbance trajectories in $\mathcal{W}$. Thus, the safety of the mixed traffic is enhanced by solving \eqref{eqn:epi_robustCount}. On the other hand, the stricter safety constraint further increases the complexity, which will be addressed in Section \ref{subsec:efficient_compu}. 
    
\begin{remark}[Uncertatinty Quantification]
    We require an accurate and non-conservative estimation of $\mathcal{W}$ for velocity error trajectories to ensure mixed traffic safety and~good~control performance. The actual disturbance trajectory should be inside or close to $\mathcal{W}$; otherwise, a gap between online prediction and real traffic behavior may still exist. A conservative estimation is not preferred either, which will shrink the feasible solution set and degrade the control performance. \hfill $\square$
\end{remark}

\begin{figure}
    \centering
    \includegraphics[width=0.48\textwidth]{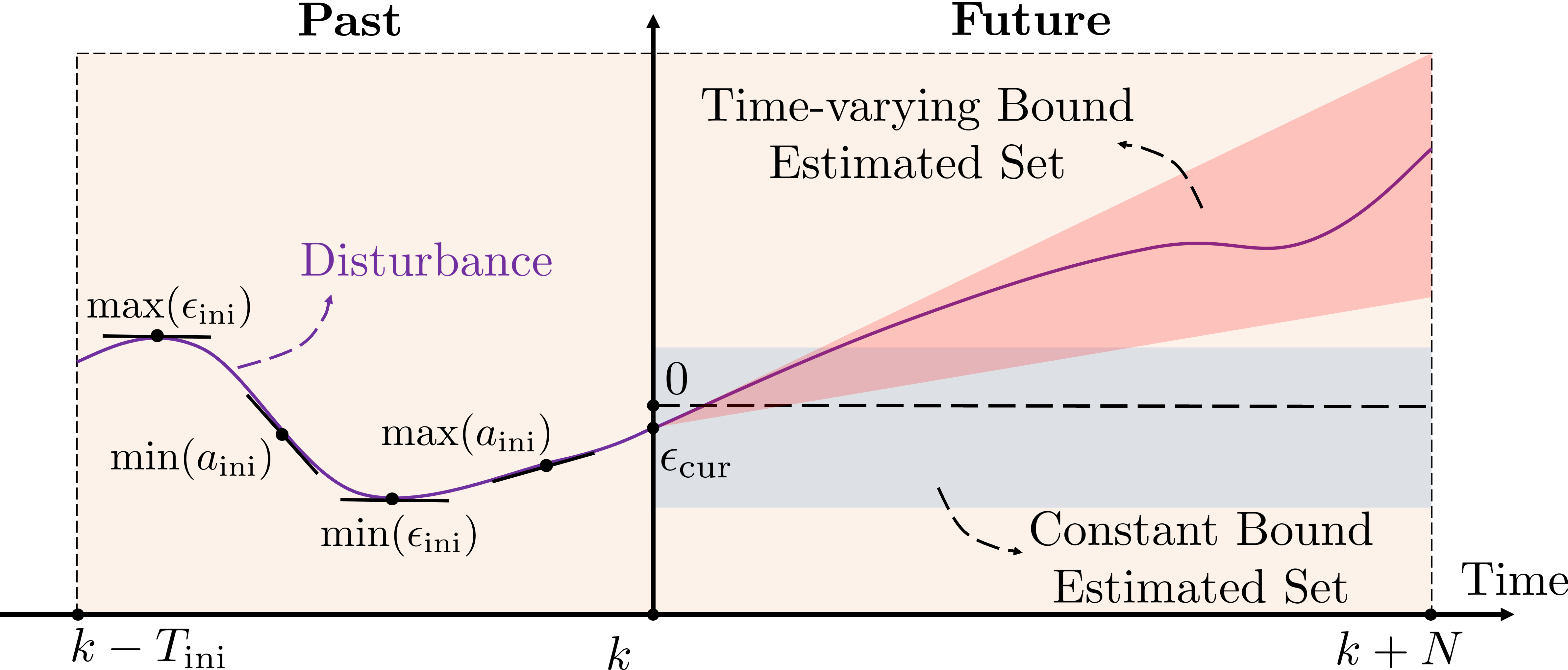}
    \caption{Schematic of two disturbance estimation methods. The purple line represents the actual disturbance trajectory and its past part is known while its future part is unknown. In the past region, the black line segment denotes the information needed for estimation. In the future region, the black dashed line represents the zero estimation while the red region and the blue region denote the time-varying bound estimated set and the constant bound estimated set, respectively.}
    \label{fig:comparison_est}
    \vspace{-2mm}
\end{figure}

\section{Disturbance Estimation and Efficient Computation}
\label{sec:Dis_Compute}
In this section, we first introduce two disturbance estimation methods based on different assumptions of human driving behaviors. We then present two solving methods of~\eqref{eqn:epi_robustCount} and compare their complexity. Also, we provide a down-sampling method of low-dimensional approximation for the disturbance set $\mathcal{W}$ for real-time computation. 

\subsection{Uncertatinty Quantification} \label{subsec:dis_est}
In our problem, the estimated disturbance set is considered an $N$-dimensional polytope that is 
\begin{equation}
\label{eq:w-subspace}
\mathcal{W} = \{\epsilon \in \mathbb{R}^N | A_{\epsilon} \epsilon \le b_{\epsilon} \},
\end{equation}
where $A_\epsilon = [I; -I]$, $b_\epsilon = [\epsilon_{\max}; -\epsilon_{\min}]$ and $\epsilon_{\max}, \epsilon_{\min}$ are the upper and lower bound vectors of $\epsilon$. The key part of estimating the disturbance set becomes estimating its (time-varying) bounds from the past velocity errors $\epsilon_{\textrm{ini}}$.

We propose two different estimation methods (see Fig.~\ref{fig:comparison_est} for illustration) and analyze their performance:

\subsubsection{Constant disturbance bounds} 
We assume that the disturbance (velocity error) of the head vehicle will not have a large deviation from its current value in a short time period based on the constant velocity model, and the disturbance variation for the future disturbance trajectory is close to its past trajectory. From the historical disturbance values $\epsilon_{\textrm{ini}}$, we can get the value of the current disturbance, i.e.,  $\epsilon_{\textrm{ini}}(\textrm{end})$, and estimate the disturbance variation as $\Delta \epsilon_{\textrm{low}} = \min(\epsilon_{\textrm{ini}})-\textrm{mean}(\epsilon_{\textrm{ini}})$ and $\Delta \epsilon_{\textrm{up}} = \max(\epsilon_{\textrm{ini}})-\textrm{mean}(\epsilon_{\textrm{ini}})$. Then, the estimated bound of the future disturbance is given~by
    \[
    \epsilon_{\min} = \epsilon_{\textrm{cur}} + \Delta \epsilon_\textrm{low}, \  
    \epsilon_{\max} = \epsilon_{\textrm{cur}} + \Delta \epsilon_\textrm{up}.
    \]

\subsubsection{Time-varying disturbance bounds} 
We can also assume the acceleration of the head vehicle will not deviate significantly from its current value based on the constant acceleration model, and its variation in the future is close to the variation in the past. We first get the past acceleration information from $\epsilon_\textrm{ini}$ as $a_\textrm{ini}(k) = \frac{\epsilon_\textrm{ini}(k+1)-\epsilon_\textrm{ini}(k)}{\Delta t}$ where $\Delta t$ is the sampling time period. Then, using a similar procedure as in the previous approach, the acceleration variation bound is estimated as $\Delta a_{\textrm{low}} = \min(a_{\textrm{ini}})-\textrm{mean}(a_{\textrm{ini}})$ and $\Delta a_{\textrm{up}} = \max(a_{\textrm{ini}})-\textrm{mean}(a_{\textrm{ini}})$. Thus, the future disturbance in an arbitrary time step $k$ is bounded by the following inequalities:
\[
\begin{aligned}
\epsilon_{\textrm{ini}}(\textrm{end}) + (a_{\textrm{cur}}+\Delta a_{\textrm{low}}) &\cdot k \Delta t \le \epsilon(k) \\ 
&\le \epsilon_{\textrm{ini}}(\textrm{end}) + (a_{\textrm{cur}}+\Delta a_{\textrm{up}}) \cdot k \Delta t.
\end{aligned}
\]

Fig. \ref{fig:comparison_est} illustrates the two disturbance estimation methods. It is clear that there exists a large gap between the actual disturbance trajectory and the zero line. For the constant disturbance bounds, the actual disturbance trajectory stays in the estimated set in the short term but will deviate from the set over time. For the second method using time-varying disturbance bounds, it includes the actual trajectory in the estimated set in this case but with a relatively conservative bound at the end of the time period. In most of our numerical simulations, the time-varying disturbance bounds outperform the constant disturbance bounds because traffic waves usually have high amplitude with low frequency.

\subsection{Efficient Computations} \label{subsec:efficient_compu}
Upon estimating $\mathcal{W}$, the robust optimization problem \eqref{eqn:epi_robustCount} is well-defined. Robust optimization is a well-studied field~\cite{lofberg2012automatic,bertsimas2011theory}. We here adapt standard robust optimization techniques to solve \eqref{eqn:epi_robustCount} and compare their complexity.


\textbf{M1: Vertex-based}.  Our first method utilizes constraints evaluated at vertices of $\mathcal{W}$ to replace the robust constraints. The compact polytope $\mathcal{W}$ can be represented as the convex hull of its extreme points as
\begin{equation}
\label{eq:w-vertex}
    \mathcal{W} = \textrm{conv}(\omega_1, \ldots, \omega_{n_\textrm{v}}),
\end{equation}
where $n_\textrm{v}$ denotes the number of extreme points, and its value is $2^N$ if no low-dimensional approximation is applied. Using this representation, we can rewrite problem \eqref{eqn:epi_robustCount} as
\begin{subequations} \label{finalformVertex}
\begin{align}
\min_{x, t} \quad& t \nonumber\\
\textrm{subject to} \quad & x_j^\tr M x_j + d^\tr x_j \le t, \,  j = 1,\ldots, n_{\textrm{v}}, \label{FFVertexC1}\\ 
& \tilde{s}_{\min} \! \le \! P_1 x_j +c_1 \! \le \! \Tilde{s}_{\max}, \, j = 1,\ldots, n_{\textrm{v}}, \label{FFVertexC2} \\
& u_{\min} \le P_2 x \le u_{\max} \label{FFVertexC3},
\end{align}
\end{subequations}
where $x_j$ represents the decision variable when the uncertainty parameter $\epsilon$ is fixed to one of the extreme points $w_j$ and the expression becomes $\textrm{col}(u, \sigma_y, w_j)$.  

\textbf{M2: Duality-based}. The second method treats robust constraint \eqref{eqn:obj_rb} the same as the first method, but forms~\eqref{eqn:safe_rb} as a sub-level optimization problem and then changes it into its dual problem to combine both levels. For example, the right hand inequality of \eqref{eqn:safe_rb} can be reformulated as 
\begin{equation}
\label{safe_opt}
\begin{aligned}
 \tilde{s}_\textnormal{max} \geq \max_{\epsilon \in \mathcal{W}} \;\; p_l ^ \tr x  + c_{1,l} , \; l = 1, \ldots, N,
\end{aligned}
\end{equation}
where $p_l^\tr$ and $c_{1,l}$ is the $l$-th row vector and element in $P_1$ and $c_{1,l}$, respectively. Given the origin representation $\mathcal{W}$ in~\eqref{eq:w-subspace}, the right-hand side of \eqref{safe_opt} is a linear program (LP). Then, we can change them to their dual problems and the strong duality of LPs ensures the new formulation is equivalent to \eqref{safe_opt}. The bi-level optimization problem becomes a min-min problem and we can combine both levels\footnote{This operation is standard; we refer the interested reader to Section 2.1 of \url{https://zhengy09.github.io/ECE285/lectures/L17.pdf}.}. The optimization problem \eqref{eqn:epi_robustCount} can then be  equivalently reformulated as 
\begin{subequations} \label{finalformDual}
\begin{align}
\min_{x_{\textrm{d}}, t, \lambda_1, \lambda_2} \quad& t \nonumber\\
\textrm{subject to} \quad 
& p_{l, \textrm{d}}^{\tr} x_{\textrm{d}}  + b_{\epsilon}^\tr \lambda_{l,1} +  c_{1,l}\le 
 \tilde{s}_{\max}, \label{FFDualC1}\\
 & A_\epsilon^\tr \lambda_{l,1} - p_{l,\epsilon} = 0,\label{FFDualC2} \\ 
&-p_{l, \textrm{d}}^\tr x_{\textrm{d}}  + b_{\epsilon}^\tr \lambda_{l,2} -  c_{1,l} \le -\tilde{s}_{\min},\label{FFDualC3} \\
 &A_\epsilon^\tr \lambda_{l,2} + p_{l,\epsilon} = 0,\label{FFDualC4} \\
&\lambda_{l,1} \geq 0, \lambda_{l,2} \geq 0, \ l = 1,2,\ldots, N,  \label{FFDualC5} \\
& \eqref{FFVertexC1}, \eqref {FFVertexC3}, \nonumber 
\end{align}
\end{subequations}
where $x_{\textrm{d}}$ is the decision variable $\textrm{col}(u, \sigma_{y})$,  and $\lambda_{l,1}, \lambda_{l,2} \! \in \! \mathbb{R}^{2n_\textrm{v}}$ are dual variables with $\lambda_1 = \textrm{col}(\lambda_{1,1}, \lambda_{2,1}, \ldots, \lambda_{N,1})$, $ \lambda_2 = \textrm{col}(\lambda_{1,2}, \lambda_{2,2}, \ldots, \lambda_{N,2})$; parameters $c_{1,l}, p_{l}$ are the same as \eqref{safe_opt} and $p_{j,\textrm{d}}$ represents $\textrm{col}(p_{l,u}, p_{l,\sigma_y})$ with  $p_l$ subdivided into $\textrm{col}(p_{l,u}, p_{l,\sigma_{y}},p_{l,\epsilon})$ corresponding to $u, \sigma_{y}$~and~$\epsilon$. 

\begin{theorem}
    Suppose \eqref{eqn:epi_robustCount} is feasible and its uncertainty set $\mathcal{W}$ is a polytope. Problems \eqref{eqn:epi_robustCount}, \eqref{finalformVertex} and \eqref{finalformDual} are~equivalent. 
\end{theorem}

The equivalence between \eqref{eqn:epi_robustCount} and \eqref{finalformVertex} is relatively straightforward. It requires standard duality arguments to establish the equivalence between \eqref{eqn:epi_robustCount} and \eqref{finalformDual}; due to the page limit, we will put the details into an extended report. 

Both \eqref{finalformVertex} and \eqref{finalformDual} are standard convex optimization problems, which can be solved using standard solvers (e.g., Mosek \cite{mosek}). We here discuss the complexity of the above two methods; see Table~\ref{Table:complexity}. The main difference lies in the different formulations of \eqref{eqn:safe_rb}, i.e.,~\eqref{FFVertexC2} and \eqref{FFDualC1} - \eqref{FFDualC5}. In \textbf{M1},~\eqref{FFVertexC2} represents $N \cdot 2^{N+1}$ inequality constraints while \eqref{FFDualC1} - \eqref{FFDualC5} together represent $2N(3N+1)$ inequality constraints in \textbf{M2}. The value $N \cdot 2^{N+1}$ is much larger than $2N(3N+2)$ when the prediction horizon $N$ is large, while there exist extra $4N^2$ decision variables in \textbf{M2}. This trade-off is also reflected in our numerical implementation. 

\begin{table}[t]
\caption{Complexity comparison between \eqref{finalformVertex} and \eqref{finalformDual}.}
\setlength{\belowcaptionskip}{0pt}
\begin{threeparttable}
    \begin{tabular}{ccc}
         \toprule
          & Decision Variables Number& Constraints Number \\
         \midrule
         \textbf{M1}  & $(n+1)T_\textrm{ini}+N+1$ & $2^{N}+N \cdot 2^{N+1}+2N$\\
         \textbf{M2} & $(n+1)T_\textrm{ini}+N+1 + 4 N^2$ & $2^{N} + 2N(3N + 2)$\\
         \midrule
         \textbf{M1 (L)}  & $(n+1)T_\textrm{ini}+N+1$ & $2^{n_\epsilon}+N \cdot 2^{n_\epsilon+1}+2N$\\
         \textbf{M2 (L)} & $(n+1)T_\textrm{ini}+N+1 + 4 N n_\epsilon$ & $2^{n_\epsilon} + 2N(3n_\epsilon + 2)$\\
         \bottomrule
         \end{tabular}
         \begin{tablenotes}
    \item[1] M1 and M2 represent the vertex-based and the duality-based method correspondingly. The last two rows denote their complexities after low-rank approximation. 
    \end{tablenotes}
    \end{threeparttable}
    \label{Table:complexity}
    \vspace{-4mm}
\end{table}
\subsection{Down-sampling strategy}
We here discuss a down-sampling strategy, adapted from~\cite{huang2021decentralized}, to relieve the exponential growth of the number of constraints. It approximates the $N$-dimensional disturbance trajectory by choosing one point for every $T_\textnormal{s}$ steps along it and performing linear interpolation. We denote the low-dimensional representation of the future disturbance trajectory as $\tilde{\epsilon} \in \mathbb{R}^{n_\epsilon}$ where $n_\epsilon = (\lfloor \frac{N-2}{T_\textnormal{s}}\rfloor+2)$. An approximated representation $\hat{\epsilon}$ of $\epsilon$ can be derived as
\[
\hat{\epsilon}^{(k)} =\left\{
\begin{aligned}
\tilde{\epsilon}^{(\bar{k}+1)}+((k-1) \ \textrm{mod} \ T_\textnormal{s}) \times 
\frac{\tilde{\epsilon}^{(\bar{k}+2)}-\tilde{\epsilon}^{(\bar{k}+1)}}{T_\textnormal{s}},\\
\quad 1 \le k \le \tilde{k}\cdot T_s \\
\tilde{\epsilon}^{(\tilde{k}+1)}+(k-\tilde{k} \cdot T_\textrm{s}-1) \times 
\frac{\tilde{\epsilon}^{(\tilde{k}+2)}-\tilde{\epsilon}^{(\tilde{k}+1)}}{N-\tilde{k} \cdot T_\textrm{s}-1},\\
 \quad \tilde{k}\cdot T_s < k \le N 
\end{aligned}\right.
\]
where $\bar{k} = \lfloor\frac{k-1}{T_\textnormal{s}}\rfloor$ and $\tilde{k} = \lfloor \frac{N-2}{T_\textrm{s}} \rfloor$. Then we can use $\tilde{\epsilon} \in \mathbb{R}^{n_\epsilon}$ to represent $\epsilon \in \mathbb{R}^N$ as 
\begin{equation}
\label{eq:low-dimen-dis}
\epsilon \approx \hat{\epsilon} = E_\epsilon \tilde{\epsilon},
\end{equation}
where $\tilde{\epsilon} \in \tilde{\mathcal{W}}$ and $\tilde{\mathcal{W}}$ can be estimated using the same methods we introduced before. Also, substituting \eqref{eq:low-dimen-dis} into our previous derivation will not affect its correctness.

The complexities of both methods \eqref{finalformVertex} and \eqref{finalformDual} after using low-dimensional approximation are updated in the last two rows of Table \ref{Table:complexity} which depend on the choice of $n_\epsilon$. Theoretically, with the same computational resource, the duality-based method allows us to choose a larger $n_\epsilon$ because the coefficient of its exponential growth term $2^{n_\epsilon}$ is $1$ while it is $2N+1$ for the vertex-based method. In our implementation, $n_\epsilon$ is usually chosen as a small number to ensure real-time computational performance and these two methods might not have obvious differences. We note that replacing $\epsilon$ with $\hat{\epsilon}$ may fail to incorporate all cases in $\mathcal{W}$ since the set of $\hat{\epsilon}$ is a subset of $\mathcal{W}$. However, our extensive simulations demonstrate that the down-sampling strategy provides satisfactory performances.

\section{Traffic Simulations}
\label{sec:Sim_Result}
In this section, we carry out nonlinear and non-deterministic traffic simulations to test the performance of robust \method{DeeP-LCC} in controlling the CF-LCC system in mixed traffic. 
Due to the page limit, we consider the time-varying bound disturbance estimation method and duality-based solving method, and the performance of other methods 
will be included in an extended report. We implemented an automatic routine transforming \eqref{finalformDual} into standard conic programs\footnote{Our open-source implementation is available at \url{https://github.com/soc-ucsd/Decentralized-DeeP-LCC/}.}, which are solved by Mosek \cite{mosek}. 

\subsection{Experimental Setup}
\begin{figure}
    \centering
    \includegraphics[width=0.48\textwidth]{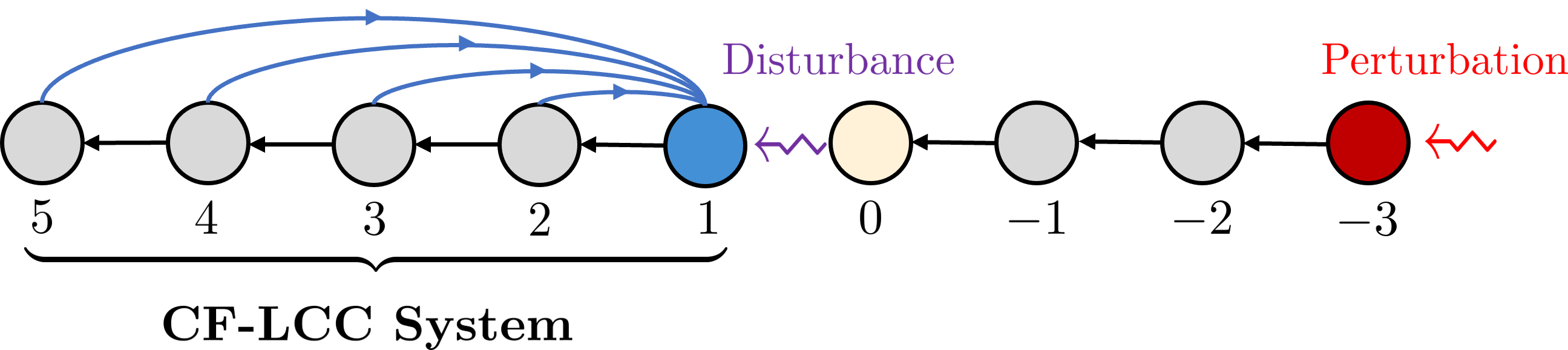}
    \setlength{\belowcaptionskip}{0pt}
    \vspace{-4mm}
    \caption{Simulation scenario. In front of the CF-LCC system, there are four preceding HDVs, where the blue node, yellow node, red node, and grey nodes represent the CAV, the head vehicle, the leading vehicle, and other HDVs, respectively.} 
    \label{fig:sim-CF-LCC}
    \vspace{-2mm}
\end{figure}

The car-following behaviors of HDVs are modeled by the nonlinear OVM model in \cite{wang2021controllability}, and a noise signal following the uniform distribution of $\mathbb{U}[-0.1,0.1]\, \textrm{m}/\textrm{s}^2$ is added to the acceleration for each HDV. For the CF-LCC system in the mixed traffic, we consider the CAV $1$ is followed by 4 HDVs, and there are three vehicles in front of the head vehicle $0$ together in the mixed traffic flow; see Fig.~\ref{fig:sim-CF-LCC} for illustration. During the simulation, a perturbation is imposed on the leading vehicle, indexed as $-3$. 

We use the following parameters in both \method{DeeP-LCC} and robust \method{DeeP-LCC}: 
\begin{enumerate}
    \item 
\textit{Offline data collection}: lengths of pre-collected data sets are $T=500$ for a small data set and $T=1500$ for a large data set with $\Delta t = 0.05 \textrm{s}$. They are collected around the equilibrium state of the system with velocity $15 \, \mathrm{m/s}$. Both $u^{\textrm{d}}$ and $\epsilon^{\textrm{d}}$ are generated by a uniform distributed signal of $\mathbb{U}[-1,1]$ which satisfies the persistent excitation requirement in Proposition \ref{proposition:data representation}; 

\item \textit{Online predictive control}: the initial signal sequence and the prediction horizon are set to $T_{\textrm{ini}} = 20$, $N=50$, respectively. For the objective function in \eqref{eqn:rDeeP}, we have $R = 0.1 I_N$ and $Q = I_N \otimes \textrm{diag}(Q_v, w_s)$ where $Q_v = \textrm{diag}(1,\ldots,1) \in \mathbb{R}^{n}$ and $w_s = 0.5$. The regularized parameters are set to $\lambda_g = 100$ and $\lambda_y = 10000$. The spacing constraints for CAV are set as $s_{\max} = 40$ m, $s_{\min}=5$ m and the bound of the spacing error is updated in each iteration as $\tilde{s}_{\max} = s_{\max}-s^*$ and $\tilde{s}_{\min} = s_{\min}-s^*$. 
\end{enumerate}
{Note that $s^*$ is also updated in each time step according to the current equilibrium state estimated by the leading vehicle's past trajectory~\cite{wang2023deep}}. The limitation of the acceleration is set as $a_{\max} = 2$ m/$\textrm{s}^2$ and $a_{\min} = -5$ m/$\textrm{s}^2$.    

\subsection{Numerical Results}
\begin{figure}[t]
\centering
\subfigure[\method{DeeP-LCC}]{\includegraphics[width=0.48\textwidth]{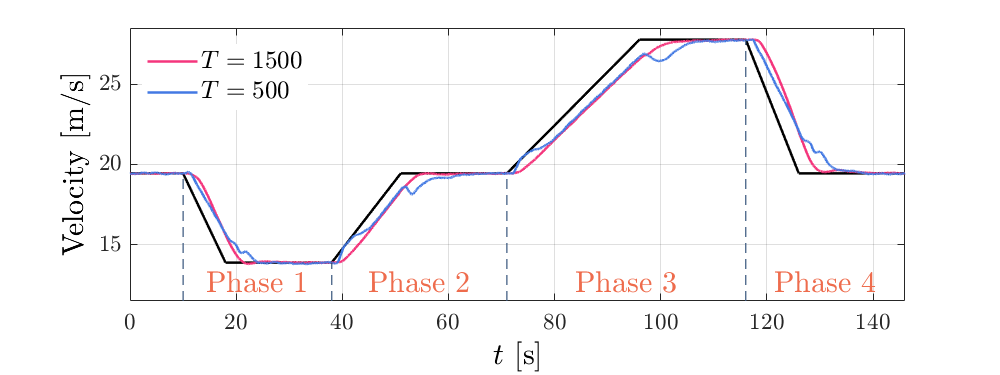}\label{subfig:Nedc_DeeP_LCC}}
\vspace{-2mm}

\subfigure[Robust \method{DeeP-LCC}]{\includegraphics[width=0.48\textwidth]{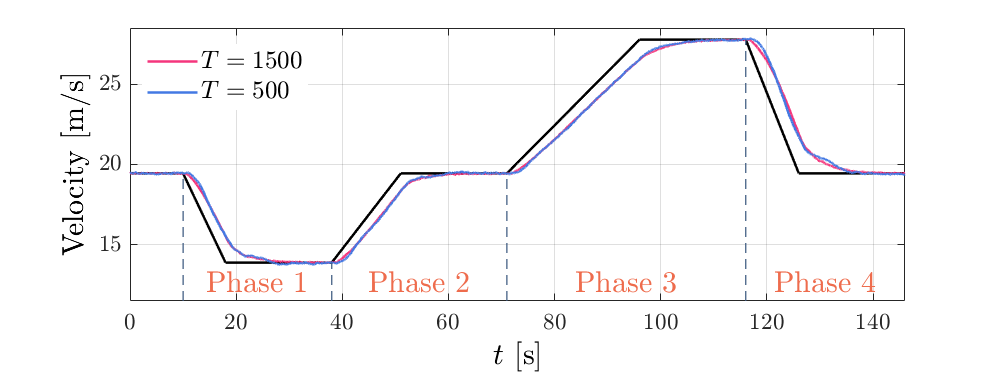}\label{subfig:Nedc_rDeeP_LCC}}
\setlength{\abovecaptionskip}{0pt}
\vspace{-1mm}
\caption{Velocity profiles in Experiment A. The black profile denotes the leading vehicle. The red profile and the blue profile represent \method{DeeP-LCC} control with data sets of size $T=1500$ and $T=500$, respectively. (a) The CAV utilizes \method{DeeP-LCC}. (b) The CAV utilizes robust \method{DeeP-LCC}.}
\label{fig:nedc_compare}
\vspace{-2mm}
\end{figure}

\textbf{Experiment A:} We first validate the control performance of robust \method{DeeP-LCC} in a comprehensive simulation scenario which is motivated by New European Driving Cycle (NEDC)~\cite{DieselNet2013}. We design the velocity trajectory of the leading vehicle as the black profile in Fig. \ref{fig:nedc_compare} and calculate the fuel consumption of the $5$ following vehicles in CF-LCC system using the numerical model in \cite{bowyer1985guide} for evaluation. 

The velocity profiles of robust \method{DeeP-LCC} and original \method{DeeP-LCC} with different sizes of data sets are shown in Fig. \ref{fig:nedc_compare}. Both methods allow for the CAV to track the desired velocity when using a large data set (see red curves in Fig. \ref{fig:nedc_compare}). However, in the case of using a small data set, the degradation of control performance for \method{DeeP-LCC} is apparent, and there are some undesired oscillations (see blue curves in Fig. \ref{subfig:Nedc_DeeP_LCC}), while robust \method{DeeP-LCC} remains a smooth velocity profile (see blue curves in Fig. \ref{subfig:Nedc_rDeeP_LCC}). Such performance degradation is highly related to the mismatch between the online prediction and real system behavior, caused by representation and estimation errors. Both original \method{DeeP-LCC} and robust \method{DeeP-LCC} employ the same data set to construct the data-driven representation \eqref{eq:CF-LCC-local-representation}, but robust \method{DeeP-LCC}  allows for a relatively small estimation error, and 
provides more margin for potential representation errors. This is one main reason that the robust \method{DeeP-LCC} performs better than \method{DeeP-LCC} for a relatively small data~set. 

Table \ref{Table:fuel_Consum} lists fuel consumption results when using the large data set. Both robust \method{DeeP-LCC} and \method{DeeP-LCC} reduce fuel consumption compared with the case with all HDVs, and the improvement in the braking phase (Phase 1 and 4) is higher than the accelerating phases (Phases 2 and 3). Moreover, we note that robust \method{DeeP-LCC} achieves better fuel economy than \method{DeeP-LCC} in all phases, $6.86\%$ vs. $3.14\%$ and $8.17\%$ vs. $4.97\%$ during Phase 1 and 4, respectively. 

\begin{table}[t]
\caption{Fuel Consumption in Experiment A (unit: $\mathrm{mL}$)}
\vspace{-2mm}
\begin{threeparttable}
    \begin{tabular}{c c c c}
         \toprule
         \noalign{\vskip-2pt}
          &All HDVs& \method{DeeP-LCC}& Robust \method{DeeP-LCC} \\
          \noalign{\vskip-2pt}
         \midrule
         Phase 1 & 145.59  &  141.02 ($\downarrow 3.14\%$) & 135.60 ($\downarrow$ \textbf{6.86\%})\\
         Phase 2 & 314.77  & 312.95 ($\downarrow 0.58\%$)  & 311.83 ($\downarrow$ \textbf{0.94\%})\\
         Phase 3 & 725.28   & 723.95 ($\downarrow 0.18\%$)  & 722.88 ($\downarrow$ \textbf{0.33\%})\\
         Phase 4 & 259.05   & 246.16 ($\downarrow 4.97\%$)  & 237.89 ($\downarrow$ \textbf{8.17\%})\\
        Total Process & 1530.15  &1509.6 ($\downarrow 1.54\%$)   & 1493.6($\downarrow$ \textbf{2.39\%})\\
         \bottomrule
    \end{tabular}
    \end{threeparttable}
    \label{Table:fuel_Consum}
    \vspace{-2mm}
\end{table}

\textbf{Experiment B:} We further validate the safety performance of robust \method{DeeP-LCC} in the braking scenario. In this experiment, the leading vehicle that moves at $15\,\mathrm{m/s}$ will brake with the maximum deceleration $-5\,\mathrm{m/s}^2$, stay at $5\,\mathrm{m/s}$ for a while, and then speed up back to $15\,\mathrm{m/s}$. 
We collect $100$ small data sets ($T=500$) and $100$ large data sets ($T=1500$) and carry out the same experiment. Recall that the safety constraint of the CAV is set from $5\,\mathrm{m}$ to $40\,\mathrm{m}$. We define ``violation" as the case where the CAV's spacing deviates more than $1\,\mathrm{m}$ from this range, and ``emergency" as the case where the spacing deviates over $5\,\mathrm{m}$ from this range. We note that, when an emergency happens, there are three possible undesired situations: 1) A rear-end collision happens; 2) The spacing of the CAV is too large which decreases the traffic capacity; 3) The controller fails to stabilize the system. 

The results are shown in Table \ref{Table:safe_rate}, which clearly shows that \method{DeeP-LCC} has a much higher violation rate and emergency rate for small data sets. Although using large data sets decreases both of them, they are still relatively high, which are $62\%$ and $51\%$ respectively. On the other hand, using the same small data sets, the robust \method{DeeP-LCC} can provide a remarkably low violation rate and emergency rate which are $5\%$ and $4\%$. Moreover, both of them are decreased to $0\%$ when using large data sets, which means perfect safety guarantees in our 100 experiments.

\begin{figure}[t]
\centering
\subfigure[Small offline data set with $T=500$]{\includegraphics[width=0.24\textwidth]{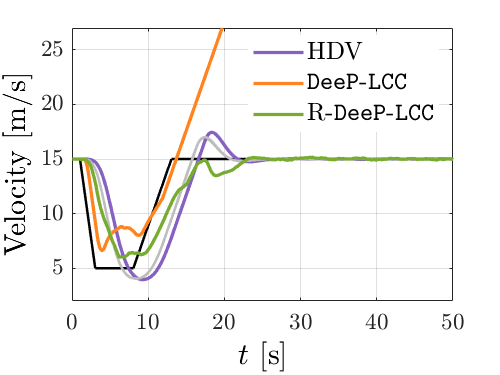}
\includegraphics[width=0.24\textwidth]{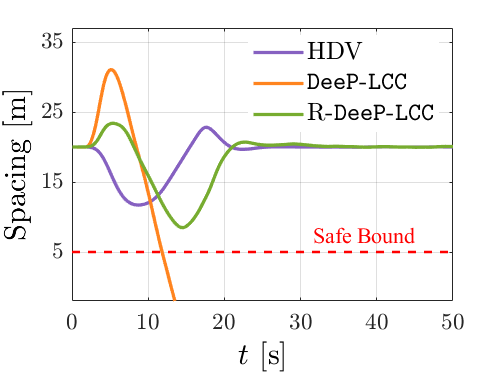}\label{fig:space_T=500}}

\vspace{-2mm}

\subfigure[Large offline data set with $T=1500$]{\includegraphics[width=0.24\textwidth]{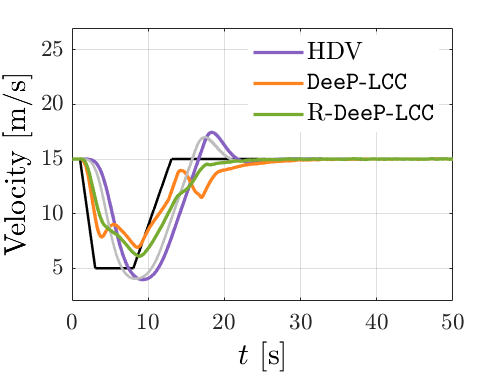}
\includegraphics[width=0.24\textwidth]{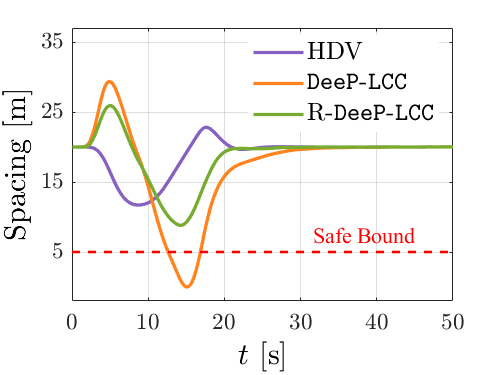}\label{fig:space_T=1500}}
\vspace{-3.5mm}
\caption{Simulation results in Experiment B. The black profile and the gray profile represent the leading vehicle and the head vehicle, respectively. The orange profile and the green profile correspond to \method{DeeP-LCC} and robust \method{DeeP-LCC}, respectively, while the purple profile corresponds to the all HDV case. (a) and (b) show the velocity and spacing profiles at different sizes of data sets. }
\label{fig:safety_compare}
\vspace{-2.5mm}
\end{figure}

\begin{table}[t]
\caption{Collision and safety Constraint Violation Rate}
\vspace{-2mm}
\begin{threeparttable}
    \begin{tabular}{c c c c c}
         \toprule
         \noalign{\vskip-2pt}
          &\multicolumn{2}{c}{\method{DeeP-LCC}}& \multicolumn{2}{c}{Robust \method{DeeP-LCC}} \\
          \cmidrule(lr{0.5em}){2-3} \cmidrule(lr{0.5em}){4-5}
          &$T = 500$ & $T = 1500$ & $T= 500$ & $T=1500$ \\
          \noalign{\vskip-2pt}
         \midrule
         Violation Rate & 74$\%$& 62$\%$& $\mathbf{5\%}$ & $\mathbf{0\%}$ \\ 
         Emergency Rate & 66$\%$& 51$\%$& $\mathbf{4\%}$ & $\mathbf{0\%}$ \\
         \bottomrule
    \end{tabular}
    \end{threeparttable}
    \label{Table:safe_rate}
    \vspace{-2mm}
\end{table}

Fig. \ref{fig:safety_compare} demonstrates two examples from small data sets and large data sets to analyze different performances of the \method{DeeP-LCC} and robust \method{DeeP-LCC}. When using a large data set, both methods exhibit smaller velocity fluctuations compared with the case of all human drivers. It can be clearly observed that the CAV controlled by robust \method{DeeP-LCC} always stays inside the safety bound for both large and small data sets, despite some small undesired velocity fluctuation for the small data set. However, \method{DeeP-LCC} is likely to lead to a rear-end collision for the small data set, and still violate the safe bound even with the large data set. Note that although the safety constraint is imposed in \method{DeeP-LCC}, it fails in the simulation due to the mismatch between the prediction and the real behavior of the system. More precisely, in prediction, \method{DeeP-LCC} considers the future velocity error of the head vehicle as $\mathbb{0}_N$ by assuming that the head vehicle accurately tracks the equilibrium velocity. Thus, the CAV decelerates or accelerates immediately when the leading vehicle starts to brake or speed up. It is, however, not the case in real-world traffic flow, and the inaccurate estimation causes the mismatch and leads to an emergency. On the other hand, robust \method{DeeP-LCC} predicts a series of the CAV's future spacing based on the estimated disturbance set and requires all of them to satisfy the safety constraint. Thus, the robust \method{DeeP-LCC} provides much stronger safety guarantees. 

\balance 
\section{Conclusion}
\label{sec:conclusion}
In this paper, we have proposed the robust \method{DeeP-LCC} for CAV control in mixed traffic. The robust formulation and disturbance set estimation methods together provide a strong safety guarantee, improve the control performance, and allow for the applicability of a smaller data set. Efficient computational methods are also provided for the real-time implementation. Extensive traffic simulations have validated the performance of robust \method{DeeP-LCC} in comprehensive and braking scenarios. Interesting future directions include learning-based estimation for future disturbances, incorporation of communication-delayed traffic data, and extension to large-scale mixed traffic scenarios.  

\bibliographystyle{IEEEtran}
\bibliography{reference}

\end{document}